\documentclass[a4paper,12pt]{amsart}
\usepackage{amssymb}
\usepackage{ifthen}
\usepackage{graphicx}
\usepackage{subfigure}
\usepackage{color,xcolor}
\nonstopmode
\numberwithin{equation}{section}
\newtheorem{thm}{Theorem}
\newtheorem{cor}{Corollary}
\newtheorem{lem}{Lemma}
\newtheorem{prop}{Proposition}

\newtheorem{conj}{Conjecture}
\newtheorem{prob}{Problem}

\theoremstyle{definition}
\newtheorem{defn}{Definition}
\newtheorem{ca}{Case}

\newtheorem{rem}{Remark}

\newenvironment{pf}[1][]{%
 \vskip 1mm
 \noindent
 \ifthenelse{\equal{#1}{}}%
  {{\slshape Proof. }}%
  {{\slshape #1.} }%
 }%
{\qed\medskip}

\newcounter{alphabet}

\newenvironment{Thm}[1][]{\refstepcounter{alphabet}%
\bigskip%
\noindent%
{\bf Theorem \Alph{alphabet}}%
\ifthenelse{\equal{#1}{}}{}{ (#1)}%
{\bf .} \itshape}{\vskip 8pt}

\newcounter{alphabet2}

\newcommand{\IN}{{\mathbb N}}

\newcommand{\ID}{{\mathbb D}}

\def\be{\begin{equation}}
\def\ee{\end{equation}}

\newcommand{\ben}{\begin{enumerate}}
\newcommand{\een}{\end{enumerate}}

\newcommand{\blem}{\begin{lem}}
\newcommand{\elem}{\end{lem}}
\newcommand{\bthm}{\begin{thm}}
\newcommand{\ethm}{\end{thm}}
\newcommand{\bcor}{\begin{cor}}
\newcommand{\ecor}{\end{cor}}
\newcommand{\beg}{\begin{exam}}
\newcommand{\eeg}{\end{exam}}
\newcommand{\begs}{\begin{examples}}
\newcommand{\eegs}{\end{examples}}
\newcommand{\bdefe}{\begin{defn}}
\newcommand{\edefe}{\end{defn}}
\newcommand{\bprob}{\begin{prob}}
\newcommand{\eprob}{\end{prob}}
\newcommand{\bques}{\begin{ques}}
\newcommand{\eques}{\end{ques}}
\newcommand{\bei}{\begin{itemize}}
\newcommand{\eei}{\end{itemize}}
\newcommand{\bcon}{\begin{conj}}
\newcommand{\econ}{\end{conj}}
\newcommand{\bop}{\begin{op}}
\newcommand{\eop}{\end{op}}

\newcommand{\bas}{\begin{assertion}}
\newcommand{\eas}{\end{assertion}}

\newcommand{\bfa}{\begin{fact}}
\newcommand{\efa}{\end{fact}}

\newcommand{\bca}{\begin{ca}}
\newcommand{\eca}{\end{ca}}

\newcommand{\bst}{\begin{step}}
\newcommand{\est}{\end{step}}

\newcommand{\bsca}{\begin{sca}}
\newcommand{\esca}{\end{sca}}

\newcommand{\bcl}{\begin{cl}}
\newcommand{\ecl}{\end{cl}}

\newcommand{\bmlem}{\begin{mlem}}
\newcommand{\emlem}{\end{mlem}}

\newcommand{\bscl}{\begin{scl}}
\newcommand{\escl}{\end{scl}}

\newcommand{\bcons}{\begin{conjs}}
\newcommand{\econs}{\end{conjs}}

\newcommand{\bprop}{\begin{prop}}
\newcommand{\eprop}{\end{prop}}

\newcommand{\br}{\begin{rem}}
\newcommand{\er}{\end{rem}}
\newcommand{\brs}{\begin{rems}}
\newcommand{\ers}{\end{rems}}
\newcommand{\bo}{\begin{obser}}
\newcommand{\eo}{\end{obser}}
\newcommand{\bos}{\begin{obsers}}
\newcommand{\eos}{\end{obsers}}
\newcommand{\bpf}{\begin{pf}}
\newcommand{\epf}{\end{pf}}
\newcommand{\ba}{\begin{array}}
\newcommand{\ea}{\end{array}}
\newcommand{\beq}{\begin{eqnarray}}
\newcommand{\beqq}{\begin{eqnarray*}}
\newcommand{\eeq}{\end{eqnarray}}
\newcommand{\eeqq}{\end{eqnarray*}}

\newcommand{\ra}{\to}

\newcommand{\ds}{\displaystyle}

\newcounter{minutes}
\divide\time by 60
\newcounter{hours}
\multiply\time by 60 \addtocounter{minutes}{-\time}

\begin{document}

\bibliographystyle{amsplain}
\title [The Bohr-type operator on analytic functions and sections ]
{The Bohr-type operator on analytic functions and sections}

\def\thefootnote{}
\footnotetext{ \texttt{\tiny File:~\jobname .tex,
          printed: \number\day-\number\month-\number\year,
          \thehours.\ifnum\theminutes<10{0}\fi\theminutes}
} \makeatletter\def\thefootnote{\@arabic\c@footnote}\makeatother

\author{Yong Huang}
 \address{Y. Huang, School of Mathematical Sciences, South China Normal University, Guangzhou, Guangdong 510631, China.}
 \email{hyong95@163.com}

\author{Ming-Sheng Liu${}^{~\mathbf{*}}$ }
 \address{M-S Liu, School of Mathematical Sciences, South China Normal University, Guangzhou, Guangdong 510631, China.} \email{liumsh65@163.com}

\author{Saminathan Ponnusamy }
\address{S. Ponnusamy, Department of Mathematics,
Indian Institute of Technology Madras, Chennai-600 036, India. }
\email{samy@iitm.ac.in}

\keywords{Bohr radius, bounded analytic functions,  Bohr-type operator, section.
  \\
${}^{\mathbf{*}}$  Corresponding author
}


\begin{abstract}
In this paper, firstly we prove two refined  Bohr-type inequalities associated with area for bounded analytic functions $f(z)=\sum_{n=0}^{\infty}a_{n}z^{n}$ in the unit disk. Later, we establish the Bohr-type operator on analytic functions and sections.
\end{abstract}

\maketitle
\pagestyle{myheadings}
\markboth{Y. Huang, M-S Liu and S. Ponnusamy}{Refined Bohr-type inequalities for bounded analytic functions}

\section{Preliminaries and some basic questions}\label{HLP-sec1}
There has been a intensive research activity on Bohr's phenomenon, examined first in 1914 by Bohr \cite{B1914}.
The main purpose of this article is to  continue the investigation on the classical Bohr inequality in the
refined formulation studied recently in \cite{IKKP2020,PW2019} for the case of analytic functions bounded in the unit disk.
See the recent survey articles \cite{AAP2016,IKKP2018,KKP2018} and \cite[Chapter 8]{GarMasRoss-2018}.
Bohr's idea naturally extends to functions of several complex variables.
The interest in the Bohr phenomena was revived in the nineties due to extensions to holomorphic functions of several complex variables and to more abstract settings.
The Bohr radius for analytic functions from the unit disk into special domains (eg. the punctured unit disk, the exterior of the closed unit disk, and concave wedge-domains) have been discussed in \cite{Abu,Abu2,Abu4,Abu3}.
Ali et al. \cite{AliBarSoly} considered Bohr's phenomenon for even and odd analytic functions and for alternating series.
This study was continued by Kayumov and Ponnusamy \cite{KP2017,KS2017}, which in turn settled one of the conjectures,  proposed in \cite{AliBarSoly}, on Bohr radius for odd analytic functions. In continuation of the investigation on this topic,
the authors in \cite{AlKayPON-19,KSS2017,LP2019} concerned the Bohr radius for the class of all sense-preserving harmonic mappings and sense-preserving  $K$-quasiconformal harmonic mappings. In \cite{PPS2002,PS2006}, authors demonstrated the classical Bohr inequality using different methods of operators. 
Several other aspects and generalizations of Bohr's inequality may be obtained from \cite{DFOOS,EPR-2017,KP2017,LPW2020,PS2004}  and the references therein for some detailed account of work on this topic. In particular, after the appearance of the articles by Abu Muhanna et al. \cite{AAP2016} and,  Kayumov and Ponnusamy \cite{KS2017}, several investigations and new problems on Bohr's inequality in the unit disk case have appeared in the literature (cf. \cite{AAL20,BhowDas-18,HLP2020, KayPon3, LP2018, LSX2018,PVW2019,PW2019}).



\subsection{Classical Inequality of H. Bohr}



Let ${\mathcal A}$ denote that class of analytic functions $f(z)=\sum_{n=0}^{\infty} a_{n} z^{n}$ in the unit disk $\mathbb{D}:=\{z \in \mathbb{C}:\,|z|<1\}$ and $\mathcal{B}=\{f\in {\mathcal A}:\,  \mbox{$|f(z)|\leq1$ in $\ID$}\}$.
For a fixed $z\in \ID$, let ${\mathcal F}_z=\{f(z):\, f\in{\mathcal A} \}$ and introduce
the Bohr operator $B_r(f)$ on ${\mathcal F}_z$, $|z| = r$,  by
$$
B_r(f):= \sum_{k=0}^{\infty}|a_k|r^k.
$$
In 1914, a classical result of Bohr \cite{B1914} states the following:

\begin{Thm}\label{Theo-A}
{\rm (Bohr \cite{B1914})}
If $f\in \mathcal{B}$ and $f(z)=\sum_{k=0}^{\infty}a_kz^k$, then
\begin{equation}
B_{r}(f) \leq 1 ~\mbox{ for all $|z|=r\leq1/3$,}
\label{liu1}
\end{equation}
and constant $1/3$ cannot be improved.
\end{Thm}

Bohr actually obtained that the inequality \eqref{liu1} is true only when $r\leq 1/{6}$. Later Riesz, Schur and Wiener, independently
established the Bohr inequality \eqref{liu1} for $r\leq 1/{3}$ and that $1/3$ is the best possible constant which is called the Bohr radius
for the space $\mathcal{B}$. 
Indeed for the function
$$\varphi_a(z)=\frac{a-z}{1-a z},\quad a\in [0,1),
$$
it follows easily that $B(\varphi_a,r)>1$ if and only if $r>1/(1+2a)$, which for $a\ra 1$ shows that $1/3$ is optimal.
Bohr's and Wiener's proofs can be found in \cite{B1914}. Other proofs of Bohr's inequality may be
found from \cite{S1927,T1962}. 
Then it is worth pointing out that there is no extremal function in $\mathcal{B}$ such that the Bohr radius
is precisely $1/3$ (cf. \cite[Corollary 8.26]{GarMasRoss-2018}).

\subsection{The Bohr inequality for bounded analytic functions}
In what follows we let $\|f\|_r^2=\sum_{n=1}^\infty|a_n|^2r^{2n}$ whenever $f(z)=\sum_{n=1}^{\infty} a_n z^n$ converges for $|z|<1$ and $r<1$. Recently, Ponnusamy  et al. \cite{PVW2019} established the following refined Bohr inequality.

\begin{Thm}\label{Theo-B}(\cite[Theorem 2]{PVW2019})
Suppose that  $f \in\mathcal{B}$, $f(z)=\sum_{n=0}^{\infty} a_{n} z^{n}$ and $f_0(z)=f(z)-a_0$. Then for $p=1,2$, the following sharp inequality holds:
$$|a_0|^p+ \sum_{n=1}^\infty |a_n|r^n+\left(\frac{1}{1+\left|a_{0}\right|}+\frac{r}{1-r}\right)\|f_0\|_r^2\leq 1 ~\mbox{for}~ r \leq \frac{1}{1+(1+|a_0|)^{2-p}}.
$$
\end{Thm}

Besides these results, there are plenty of works about Bohr's phenomenon. As a consequence of the development on this topic, the notion of Rogosinski's inequality and Rogosinski's radius were investigated in  \cite{EDB1986,WR1923,IG1925}. Furthermore, Kayumov and Ponnusamy \cite{KP2019} introduced and studied the Bohr-Rogosinski inequality and found the Bohr-Rogosinski radius. Based on the initiation of \cite{KayPon3}, several forms of Bohr-type inequalities for the family $\mathcal{B}$
were considered in \cite{LSX2018} when the Taylor coefficients of classical Bohr inequality are partly or completely replaced by higher order derivatives of $f$. Let $S_r(f)$ denote the area of the image of the subdisk $|z|<r$ under the mapping $f$ and
when there is no confusion, we let for brevity $S_r$ for $S_r(f)$. Let us now recall a couple of recent results for our reference.

\begin{Thm}\label{Theo-C}
{\rm (\cite{LLP2020})}
Suppose that  $f \in\mathcal{B}$, $f(z)=\sum_{n=0}^{\infty} a_{n} z^{n}$ and $f_0(z)=f(z)-a_0$. Then for $p=1,2$, the following sharp inequality holds:
$$|a_0|^p+ \sum_{n=1}^\infty |a_n|r^n +\left(\frac{1}{1+\left|a_{0}\right|}+\frac{r}{1-r}\right)\|f_0\|_r^2+
\left (\frac{8}{9}\right )^{3-2p}\frac{S_{r}}{\pi} \leq 1 ~\mbox{for}~ r \leq \frac{1}{2+(1-|a_0|)^{p-1}}.
$$
\end{Thm}

Refined version of Theorem~A, 
without the third term as in Theorem~C 
were obtained in \cite{IKKP2020}. See also \cite[Remarks 1, 2 and 3]{IKKP2020} and  \cite[Theorem 2]{IKKP2020} for the analog of this theorem
with a replacement of the constant term by $|f(z)|^2$.

\begin{Thm}\label{VINITI-1}
Suppose that  $f \in\mathcal{B}$,  $f(z)=\sum_{n=0}^{\infty} a_{n} z^{n}$,
and  $S_r$ denotes the area of the image of the subdisk $|z|<r$ under the mapping $f$.
Then
$$\sum_{k=0}^\infty |a_k|r^k+\frac{16}{9}\left (\frac{S_r}{\pi}\right )+\lambda\left (\frac{S_r}{\pi}\right )^2  \leq 1 ~\mbox{ for  }~ r \leq \frac{1}{3}
$$
where
 $$
 \lambda= \frac{4 (486 - 261 a - 324 a^2 + 2 a^3 + 30 a^4 + 3 a^5)}{ 81 (1 + a)^3 (3 - 5 a)}=18.6095 \ldots
 $$ and  $a\approx 0.567284$, is the unique positive root of the equation $\psi (t)=0$ in the interval $(0,1)$, where
$$ \psi (t)= -405 + 473 t + 402 t^2 + 38 t^3 + 3 t^4 + t^5.
$$
The equality is achieved for the function
$$f(z) = \frac{a-z}{1-az}.
$$
\end{Thm}

In view of Theorems~B and C, it is natural to raise the following.

\bprob\label{HLP-prob1}
Whether we can derive sharp version of Theorem~C 
in the setting of Theorem~D 
with the additional non-negative term and without decreasing the radius?
\eprob

In Theorems \ref{HLP-th1} and \ref{HLP-th2}, we present an affirmative answer to this question.

\subsection{The Bohr operator on analytic functions and sections}
The notion closely related to the Bohr radius is the Rogosinski radius contained in the following result of Rogosinski \cite{WR1923}.

\begin{Thm}\label{Theo-D}(\cite{WR1923})
If $f(z)=\sum_{n=0}^{\infty} a_{n} z^{n}\in \mathcal{B} $ and $|f(z)|\leq 1$ for all $|z|\leq r$, then for every
$k\in \IN_{0}=\{0,1,2,\cdots\}$ and $0<r\leq 1$, each section $s_{k}(f):=s_{k}(z; f)=\sum_{n=0}^{k} a_{n} z^{n}$ of $f$ satisfies
the inequality
$$
\left|s_{k}(f)\right| \leq 1
$$
for $|z|\leq r/2 $. The constant $r/2$ cannot be improved.
\end{Thm}

The number $r/2$ (stated actually in \cite{WR1923} with $r=1$) in Theorem~E 
is known as the Rogosinski radius. In \cite{AES2005,A2012}, Aizenberg et al. extended the Rogosinski phenomenon to holomorphic functions of 
several complex variables. Aizenberg \cite{A2012} also studied the Rogosinski radius on Hardy spaces and  Reinhardt domains. For more recent advances, 
see for example \cite{AlKayPON-19,BK1997,B1914}.


In \cite{BhowDas-18}, Bhowmik and Das established the following result: {\it if $f\prec h$ in $\mathbb{D}$, then
\be \label{LHS-eq1}
B_{r}(f)\leq B_{r}(h),~~0\leq r\leq \frac{1}{3}.
\ee
}
Later in \cite{AlKayPON-19} this inequality has been extended for the family of quasi-subordinations which includes both
subordination and majorization. Alternate proof of \eqref{LHS-eq1} has been obtained in \cite{AAL20}. More recently,
Gang and Ponnusamy  refined this result and in particular, they derived the following sharper version of \eqref{LHS-eq1}:

\begin{Thm}\label{Theo-E} {\rm (\cite[Corollary 2]{LP2020})}
Suppose that $f\prec h$ in $\mathbb{D}$ and
\be \label{LHS-eq2}
r_1(x)=
\begin{cases}
\displaystyle \sqrt{\frac{1-x}{2}} & \text{for}\quad x\in[0,\frac{1}{2}),\\[4mm]
 \displaystyle \frac{1}{1+2x} & \text{for}\quad x\in[\frac{1}{2},1].
\end{cases}
\ee
Then we have
\begin{enumerate}
\item[{\rm (a)}] $\ds  B_{r}(f)\leq B_{r}(h)  $ $\ds
\text{ for } r \leq r_1(|f'(0)/h'(0)|),$
when $h'(0) \neq0$.
\item[{\rm (b)}] $\ds  B_{r}(f)\leq B_{r}(h)  $ $\ds
\text{ for } r \leq 1/3,$
when $h'(0)=0$.
\end{enumerate}
Moreover, $r_1(|f'(0)/h'(0)|)$ cannot be improved if  $|f'(0)/h'(0)|\in[1/2,1)\cup\{0\}$, and the constant $1/3$
 in {\rm (b)} cannot be improved.
\end{Thm}

Note that $r_1(x)\ge 1/3$ for $x\in [0,1]$.

%

The paper is organized as follows. In Section \ref{HLP-sec2}, we present the main results of this paper. In
Theorems \ref{HLP-th1} and \ref{HLP-th2}, we present an affirmative answer to Problem \ref{HLP-prob1}. Also,
we establish the Bohr-type operator on analytic functions and also on their sections. In Section \ref{HLP-sec3},
we state and prove several lemmas which are needed for the proofs of the following theorems.
In Sections \ref{HLP-sec4} and \ref{HLP-sec5}, we present the proofs of the main results.

\setcounter{equation}{0}
\section{Main Results}\label{HLP-sec2}


\bthm\label{HLP-th1}
Suppose that  $f \in\mathcal{B}$,  $f(z)=\sum_{n=0}^{\infty} a_{n} z^{n}$, $f_0(z)=f(z)-a_0$, and  $S_r$ denotes the area of the image of the subdisk $|z|<r$ under the mapping $f$.
Then
$$
A(r):=\sum_{n=0}^{\infty}\left|a_{n}\right| r^{n}+\left(\frac{1}{1+\left|a_{0}\right|}+\frac{r}{1-r}\right)  \|f_0\|_r^2+\frac{8}{9}\left(\frac{S_{r}}{\pi}\right)+\lambda\left(\frac{S_{r}}{\pi}\right)^{2} \leq 1 ~for~ r \leq \frac{1}{3}
$$
where
$$
\lambda=\frac{-2673+2502 a_{*}+2025 a_{*}^{2}-332 a_{*}^{3}-255 a_{*}^{4}+6 a_{*}^{5}+7 a_{*}^{6}}{162(a_{*}+1)^{3}(5 a_{*}-3)}\approx 14.796883 ,
$$
and $a_{*}\approx 0.587459$, is the unique root in $(0, 1)$ of the equation $\psi(t)=0$, where
\begin{eqnarray}\label{liu21}
\psi(t)&=&2 t^{8}+4 t^{7}-4 t^{6}-76 t^{5}-1344 t^{4}-1396 t^{3}\\
&&+12420 t^{2}+15804 t-13122.\nonumber
\end{eqnarray}
The equality in the last inequality is achieved for the function $f(z)=\frac{a-z}{1-a z}$.
\ethm

\bthm \label{HLP-th2}
Assume the hypotheses of Theorem \ref{HLP-th1}. Then
$$
B(r):=\left|a_{0}\right|^{2}+\sum_{n=1}^{\infty}\left|a_{n}\right| r^{n}+\left(\frac{1}{1+\left|a_{0}\right|}+\frac{r}{1-r}\right) \|f_0\|_r^2+\frac{9}{8}\left(\frac{S_{r}}{\pi}\right)+\mu\left(\frac{S_{r}}{\pi}\right)^{2} \leq 1
$$
for $r \leq \frac{1}{3-\left|a_{0}\right|}$, where
$$
\mu=\frac{-80919+119556 a_{**}-57591 a_{**}^{2}+11664 a_{**}^{3}-1620 a_{**}^{4}}{8(a_{**}+1)^{2}(a_{**}-3)^{3}\left(9 a_{**}^{3}-33 a_{**}^{2}+29 a_{**}-1\right)}=13.966088\cdots ,
$$
and $a_{**}\approx 0.638302$ is the unique root in $(0, 1)$ of the equation $\phi(t)=0$, where
\begin{eqnarray}\label{liu22}
\phi(t)&=&-1296 t^{8}+17172 t^{7}-154386 t^{6}+798660 t^{5}-2361960 t^{4}\\
&&+4132944 t^{3}-4244238 t^{2}+2344464 t-524880.\nonumber
\end{eqnarray}
The equality in the last inequality is achieved for the function $f(z)=\frac{a-z}{1-a z}$.
\ethm

In the next three theorems, we use the idea of \cite{AAL20} to establish some similar results for derivatives.

\bthm \label{HLP-th3}
Suppose $f\prec h$ in $\mathbb{D}$.
Then for each $k$-th sections of $f'$ and $h'$, we have the following
\begin{enumerate}
\item[{\rm (1)}] $\ds
\left|s_{k}\left(f^{\prime}\right)\right| \leq  B_{r}\left(s_{k}\left(h^{\prime}\right)\right)$, for  $0 \leq r \leq \frac{1}{2}\Big(1-\sqrt{\frac{2}{3}}\Big)$
\item[{\rm (2)}] $\ds  B_{r}\left(s_{k}\left(f^{\prime}\right)\right) \leq  B_{r}\left(s_{k}\left(h^{\prime}\right)\right)$ for $0 \leq r \leq \frac{1}{3}\Big(1-\sqrt{\frac{2}{3}}\Big)$.
\end{enumerate}
\ethm

%


\bthm \label{HLP-th6}
Suppose $h$ is analytic in $\mathbb{D}$. If $0 \leq r \leq \frac{1}{3}\Big(1-\sqrt{\frac{2}{3}}\Big)$, then
$$
 B_{r}\left((h \circ \varphi)^{\prime}\right) \leq  B_{r}\left(h^{\prime}\right)
$$
for every Schwarz function $\varphi$.
\ethm

\bthm \label{HLP-th7}
Suppose the analytic functions $f,g$ and $h$ satisfy $f(z)=g(z) h(\varphi(z)), z\in \mathbb{D}$, for some Schwarz function $\varphi$. Further suppose that $|g(z)| \leq b$ for $|z| <1$. Then
$$
 B_{r}\left({f^{\prime}}\right) \leq b\left(2 B_{r}(h)+ B_{r}\left(h^{\prime}\right)\right)
$$
holds for $|z|=r \leq r_{2}$, where $r_{2}=1-\sqrt{\frac{2}{3}}$.
\ethm

\section{Key lemmas and their Proofs}\label{HLP-sec3}
In order to establish our main results, we need the following lemmas. 


\blem\label{HLP-lem1} (\cite{LLP2020})
Suppose that  $f \in\mathcal{B}$ and $f(z)=\sum_{n=0}^{\infty} a_{n} z^{n}$. Then for any $ N\in \mathbb{N}$, the following inequality holds:
$$
\sum_{n=N}^{\infty}\left|a_{n}\right| r^{n}+{sgn}(t) \sum_{n=1}^{t}\left|a_{n}\right|^{2} \frac{r^{N}}{1-r}+\left(\frac{1}{1+\left|a_{0}\right|}+\frac{r}{1-r}\right) \sum_{n=t+1}^{\infty}\left|a_{n}\right|^{2} r^{2 n} \leq\left(1-\left|a_{0}\right|^{2}\right) \frac{r^{N}}{1-r}
$$
for $r\in[0,1)$, where $t=[(N-1) / 2]$.
\elem

\blem\label{HLP-lem2}
{\rm (\cite[Lemma 1]{KP2017})}
Suppose that  $f \in\mathcal{B}$,  $f(z)=\sum_{n=0}^{\infty} a_{n} z^{n}$,
and  $S_r$ denotes the area of the image of the subdisk $|z|<r$ under the mapping $f$.
Then the following sharp inequality holds:
$$
\frac{S_{r}(f)}{\pi}:=\sum_{n=1}^{\infty} n\left|a_{n}\right|^{2} r^{2 n} \leq r^{2} \frac{\left(1-\left|a_{0}\right|^{2}\right)^{2}}{\left(1-\left|a_{0}\right|^{2} r^{2}\right)^{2}} \text { for } 0<r \leq 1 / \sqrt{2}.
$$
\elem

In \cite{BD2019}, Bhowmik et al. established the Bohr phenomenon for the derivative of an analytic self map $\varphi$ of $\mathbb{D}$,
where $\varphi(0)=0$. The next result is stated with radius dependence on the initial coefficient of $\varphi$ and the proof is on the
lines of the proof of \cite{BD2019}.

\blem\label{Theo-F}
Let $\varphi\in {\mathcal B}$ such that $\varphi (0)=0$, i.e., $\varphi $ is a Schwarz function.
Then
$$ B_{r}\left(\varphi ^{\prime}\right) \leq 1 \mbox{ for $|z|= r \leq r_{0}(|\varphi '(0)|)$}, ~r_0(x)=1-\sqrt{\frac{1+x}{2+x}}.
$$
Here $r_0(x)\geq r_0(1)=1-\sqrt{\frac{2}{3}}$ for all $x\in [0,1]$ and the number $r_{0}(1)$ is optimal.
\elem

\bpf
We may write $\varphi(z)=z\omega (z)$, where  $\omega$ is an analytic self map of $\mathbb{D}$ with $\omega (0)=\varphi'(0)$.
Let $\omega (z)=\sum_{n=0}^{\infty} w_{n} z^{n}$. Then $|w_{n}|\leq 1-|w_0|^2$ for $n\geq 1$ and
$\varphi^{\prime}(z)=z\omega'(z) +\omega (z)
$
so that
\beqq
 B_{r}\left(\varphi ^{\prime}\right)
&=& |w_0|+\sum_{n=1}^{\infty} (n+1)|w_{n}| r^{n}\\
&\leq & |w_0|+(1-|w_0|^2)\left [\frac{1}{(1-r)^2}-1\right ] \\
&=& 1-(1-|w_0|^2)\left [\frac{2+|w_0|}{1+|w_0|}-\frac{1}{(1-r)^2}\right ],
\eeqq
which is less than or equal to $1$ provided $r\leq  r_{0}(|w_0|)=r_{0}(|\varphi '(0)|)$. It is a simple exercise to see that
for $\zeta (z)=z\left (\frac{z-a}{1-az}\right )$, $0\leq a<1$, we have
$$
B_{r}\left(\zeta ^{\prime}\right) =a+ (1-a^2)\frac{2r-ar^2}{(1-ar)^2}>1
~\mbox{ whenever }~r>\frac{1}{a}\left [ 1 - \sqrt{\frac{1+a}{1+2a}}\right ].
$$
Allowing $a\rightarrow 1$, we see that $r_0(1)=1-\sqrt{\frac{2}{3}}$ is optimal.
\epf

\br
The function $r_0(x)\, (x\in [0,1)$) is strictly decreasing from $1-\sqrt{\frac{1}{2}}$ to $1-\sqrt{\frac{2}{3}}$.
\er

\blem\label{HLP-lem3}
If $\varphi$ is a Schwarz function, then $ B_{r}(\varphi) \leq 1$ for $|z|=r\leq r_1(|\varphi '(0)|)$, where
$r_1(x)$ is defined by \eqref{LHS-eq2} and $r_1(|\varphi '(0)|)$ cannot be improved if $|\varphi'(0)|\in[1/2,1)\cup\{0\}$.
\elem

\bpf
We may write $\varphi(z)=z\omega (z)$, where  $\omega\in {\mathcal B}$ and $\omega (0)=\varphi'(0)$. According to \cite[Theorem 1]{LP2020},
$$ B_{r}(\omega) \leq 1 ~\mbox{ for $|z|\leq r_1(|\omega (0)|)=r_1(|\varphi '(0)|)$},
$$
and $r_1(x)$ is as in the statement. Note that $|\varphi '(0)|\in [0,1]$ and the function $r_1(x)$
$(x\in [0,1)$) is strictly decreasing from $1/\sqrt{2}$ to $1/3$. The sharpness part follows as
in \cite[Theorem 1]{LP2020}.
\epf

Clearly, Lemma \ref{HLP-lem3}  refines \cite[Lemma 1]{AAL20}.

\blem\label{HLP-lem4}
Let $\varphi$ be a Schwarz function, $j\in \mathbb{N}$,  $r_{0}(x)$ and $r_{1}(x)$ be as in Lemma \ref{Theo-F} and \eqref{LHS-eq2}, respectively.
Then the following inequalities hold:
\be\label{LHS-eq3}
\sup _{|z|\leq r}\left|s_{k}\left(\varphi^{\prime}(z) \varphi^{j}(z)\right)\right| \leq|z|^{j}
~\mbox{ for }~r\leq \frac{1}{2}r_{0}(|\varphi '(0)|) 
\ee
and
\be\label{LHS-eq4}
 B_{r}\left(s_{k}\left(\varphi^{\prime}(z) \varphi^{j}(z)\right)\right) \leq|z|^{j} ~\mbox{ for }~r\leq r_{0}(|\varphi '(0)|)  r_{1}(|\varphi '(0)|^j). 
\ee
\elem
\bpf
First we remark once again that  $r_{0}(|\varphi '(0|)\geq r_{0}(1)=1-\sqrt{\frac{2}{3}} $ and  $ r_{1}(|\varphi '(0)|^j) \geq r_{1}(1)=1/3$.
Next, we observe (by Lemma \ref{Theo-F}) that
$$
\left|\varphi^{\prime}(\varphi / z)^{j}\right|=\left|\varphi^{\prime} |\, | \varphi / z\right|^{j} \leq  B_{r}\left(\varphi^{\prime}\right)|\varphi / z|^{j} \leq|\varphi / z|^{j} \leq 1
$$
for $|z|= r \leq r_{0}(|\varphi '(0)|)$. According to Theorem~E, 
we obtain \eqref{LHS-eq3}.
Again
$$ B_{r}\left(\varphi^{\prime}(\varphi / z)^{j}\right) \leq  B_{r}\left(\varphi^{\prime}\right)
 B_{r}\left( (\varphi / z)^{j}\right) \leq  B_{r}\left( (\varphi / z)^{j}\right) ~\mbox{ for }~ 0 \leq r\leq r_{0}(|\varphi '(0)|).
$$
which in turn implies  (cf.  \cite[Theorem 1]{LP2020}) that
$$ B_{r}\left(\varphi^{\prime}(\varphi / z)^{j}\right) \leq  1 ~\mbox{ for }~  0 \leq r\leq r_{0}(|\varphi '(0)|)  r_{1}(|\varphi '(0)|^j).
$$
As
$ B_{r}\left(s_{k}\left(\varphi^{\prime}(\varphi / z)^{j}\right)\right) \leq  B_{r}\left(\varphi^{\prime}(\varphi / z)^{j}\right)$, the desired conclusion follows from the last inequality.
\epf

\blem\label{HLP-lem5}
There is a unique root $a_{*}\approx 0.587459$ in $(0, 1)$ of the equation $ \psi(a)=0$, where $\psi$ is given by \eqref{liu21}.
\elem

\bpf
Firstly, we begin by a direct computation that $\psi(0)=-13122<0$ and $\psi(1)=12288>0$. Also,
\begin{eqnarray*}
\psi^{\prime}(a)&=&16 a^{7}+28 a^{6}-24 a^{5}-380 a^{4}-5376 a^{3}-4188 a^{2}+24840 a+15804,\\
\psi^{\prime \prime}(a)&=&112 a^{6}+168 a^{5}-120 a^{4}-1520 a^{3}-16128 a^{2}-8376 a+24840,\\
\psi^{(3)}(a)&=&672 a^{5}+840 a^{4}-480 a^{3}-4560 a^{2}-32256 a-8376,\\
\psi^{(4)}(a)&=&3360 a^{4}+3360 a^{3}-1440 a^{2}-9120 a-32256,\\
\psi^{(5)}(a)&=&13440 a^{3}+10080 a^{2}-2880 a-9120, ~\mbox{ and,}\\
\psi^{(6)}(a)&=&40320 a^{2}+20160 a- 2880=40320 \Big(a^{2}+\frac{1}{2} a-\frac{1}{14}\Big)\\
&=& 40320 \Big(a+\frac{\sqrt{105}+7}{28}\Big)\Big(a-\frac{\sqrt{105}-7}{28}\Big).
\end{eqnarray*}

Since 
$\psi^{(6)}(a)\leq 0$ holds if $0\leq a<a_{1} $, where $a_{1}=\frac{\sqrt{105}-7}{28}\approx 0.115963$, and
$\psi^{(6)}(a)\geq 0$ holds if $a_{1}\leq a\leq 1$, we obtain that $\psi^{(5)}(a)$ is a decreasing function of
$a\in [0,a_{1}]$ and is an increasing function of $a\in [a_{1},1]$.

Note that $\psi^{(5)}(0)=-9120<0$ and $\psi^{(5)}(1)=11520>0$, we obtain that $\psi^{(5)}(a)<0$ for $[0, a_1]$, and $\psi^{(5)}(a)$ has a unique zero $a_2\in (a_1, 1)$. Thus, 
we find that $\psi^{(5)}(a)<0$ for $0<a<a_{2} $ and $\psi^{(5)}(a)>0$  for $a_{2}<a<1$. 
Therefore $\psi^{(4)}(a)$ is a decreasing function of $a\in [0,a_{2})$ and is an increasing function of
$a\in (a_{2},1]$. Since
$$
\psi^{(4)}(0)=-32256<0 ~\mbox{ and }~
\psi^{(4)}(1)=-36096<0,
$$
we obtain that $~\psi^{(4)}(a)<0$ for $a\in [0, 1]$, which shows that $\psi^{(3)}(a)$ is decreasing in $[0, 1]$ so that
$$\psi^{(3)}(a)\leq \psi^{(3)}(0)=-8376<0 ~\mbox{  for $a\in [0, 1]$}.
$$
This implies that $\psi ''(a)$ is a monotonically decreasing function of $a$ for $a\in [0, 1]$.

Note that $\psi''(0)=24840>0$ and $\psi ''(1)=-1024<0$, we obtain that $\psi{''}(a)$ has a unique zero $a_3\in (0, 1)$. 
Thus, $\psi^{\prime}(a)$ is an increasing function of $a\in (0,a_{3})$ and is a decreasing function of $a\in (a_{3},1)$. Since
$$
\psi^{\prime}(0)=15804>0 ~\mbox{ and }~ \psi^{\prime}(1)=30720>0,
$$
it follows that $\psi^{\prime}(a)>0$ for $a\in [0, 1]$. We conclude that $\psi(a)$ is an increasing function of $a$ in $[0, 1]$,
with  $\psi(0)=-13122<0$ and $\psi(1)=12288>0$ which yields that the equation $ \psi(a)=0$ has a unique positive root in $(0,1)$.
Through the calculation of Mathematica or Maple, we obtain that the root is $a^{*}\approx 0.587459$.
\epf

\blem\label{HLP-lem6}
There is a unique root $a_{**}\approx 0.638302$ in $(0, 1)$ of the equation $ \phi(a)=0$, where $\phi $ is given by \eqref{liu22}.
\elem

\bpf
Firstly  we observe that $\phi(0)=-524880<0$ and $\phi(1)=6480>0$. By differentiating $ \phi(a)$, we easily obtain
\begin{eqnarray*}
\phi^{\prime}(a)&=&-10368 a^{7}+120204 a^{6}-926316 a^{5}+3993300 a^{4}-9447840 a^{3}\\
&&+12398832 a^{2}-8488476 a+2344464,\\
\phi^{\prime \prime}(a)&=&-72576 a^{6}+721224 a^{5}-4631580 a^{4}+15973200 a^{3}-28343520 a^{2}\\
&&+24797664 a-8488476,\\
\phi^{(3)}(a)&=&-435456 a^{5}+3606120 a^{4}-18526320 a^{3}+47919600 a^{2}-56687040 a\\
&&+24797664,\\
\phi^{(4)}(a)&=&-2177280 a^{4}+14424480 a^{3}-55578960 a^{2}+95839200 a-56687040,\\
\phi^{(5)}(a)&=&-8709120 a^{3}+43273440 a^{2}-111157920 a+95839200, \\
\phi^{(6)}(a)&=&-26127360 a^{2}+86546880 a-111157920.
\end{eqnarray*}

Since the discriminant of the equation $\phi^{(6)}(a)=0$ is less than $0$ and the coefficient of $a^2$ in
$\phi^{(6)}(a)$ is negative, we obtain that $\phi^{(6)}(a)<0$
for $a\in [0, 1]$. Thus, $\phi^{(5)}(a)$ is a decreasing function of $a$ for $a\in [0, 1]$, and therefore,
$$\phi^{(5)}(a)\geq \phi^{(5)}(1)=19245600>0 ~\mbox{ for $a\in [0, 1]$,}
$$
showing that $\phi^{(4)}(a)$ is increasing in $[0, 1]$ so that
$$\phi^{(4)}(a)\leq \phi^{(4)}(1)=-4179600<0 ~\mbox{  for $a\in [0, 1]$}.
$$
This implies that $\phi^{(3)}(a)$ is a decreasing function of $a$ for $a\in [0, 1]$, so that
$$\phi^{(3)}(a)\geq \phi^{(3)}(1)=674568>0 ~\mbox{ for $a\in [0, 1]$.}
$$
We conclude that $\phi^{\prime\prime}(a)$ is an increasing function of $a$ in $[0, 1]$ and therefore,
$$\phi^{\prime \prime}(a) \leq \phi^{\prime \prime}(1)=-44064<0 ~\mbox{ for $a\in [0, 1]$,}
$$
which shows that $\phi^{\prime}(a)$ is decreasing in the interval $[0, 1]$.
As $\phi^{\prime}(0)>0$ and $\phi^{\prime}(1)=-16200<0$, there exists a unique $a_4\in (0,1)$ (in fact $a_4\approx 0.853918$)
such that $\phi^{\prime}(a)>0$ when $0<a<a_4$, and $\phi^{\prime}(a)<0$ when $a_4<a<1$. Thus, $\phi(a)$ is an increasing
function of $a\in (0,a_4)$ and is a decreasing function of $a\in (a_4, 1)$ with $\phi (0)<0$ and $\phi (1)=6480>0$.
Thus, the equation $ \phi(a)=0$ has a unique positive root $a_{**}$ in $(0,1)$. We find that $a_{**}\approx 0.638302$ as desired.
\epf

\section{Bohr-type inequalities for bounded analytic functions}\label{HLP-sec4}

\subsection{Proof of Theorem \ref{HLP-th1}}
Let $|a_{0}|=a \in(0,1)$. By Lemma \ref{HLP-lem1} with $N=1$, and  Lemma \ref{HLP-lem2}, we have
$$
A(r) \leq a+\frac{\left(1-a^{2}\right) r}{1-r}+\frac{8\left(1-a^{2}\right)^{2} r^{2}}{9\left(1-a^{2} r^{2}\right)^{2}}+\lambda \frac{\left(1-a^{2}\right)^{4} r^{4}}{\left(1-a^{2} r^{2}\right)^{4}}:=A_{1}(r)
$$
Since $A_{1}(r)$ is an increasing function of $r$, we have for $r\leq 1/3$,
\begin{eqnarray*}
A(r) &\leq& A_{1}\left(\frac{1}{3}\right)=a+\frac{1-a^{2}}{2}+\frac{8\left(1-a^{2}\right)^{2}}{\left(9-a^{2}\right)^{2}}+81 \lambda \frac{\left(1-a^{2}\right)^{4}}{\left(9-a^{2}\right)^{4}}\\
&=&1-\frac{(1-a)^{3}}{2\left(9-a^{2}\right)^{4}} A_{2}(a),
\end{eqnarray*}
where
\begin{eqnarray*}
A_{2}(a)&=&5265+2673 a-1251 a^{2}-675 a^{3}+83 a^{4}+51 a^{5}-a^{6}-a^{7}\\
&&+162 \lambda\left(-1-3 a-2 a^{2}+2 a^{3}+3 a^{4}+a^{5}\right),\\
&=&5265+2673 a-1251 a^{2}-675 a^{3}+83 a^{4}+51 a^{5}-a^{6}-a^{7}\\
&& -162 \lambda (1-a^{2})(1+a)^{3}.
\end{eqnarray*}

Next, we will verify that the function $A_{2}(a)$ has exactly one stationary point $a_{*}=0.587459\cdots$ in $[0, 1]$, which is the unique root in $[0,1]$ of the equation $ A_{2}(a)=0$. We remark that $ A_{2}'(a)=0$ is equivalent to
$$
162\lambda (a+1)^{3}(5 a-3)=-2673+2502 a+2025 a^{2}-332 a^{3}-255 a^{4}+6 a^{5}+7 a^{6}  ,
$$
from which we obtain the value of $\lambda$ mentioned in the statement of the theorem.

In fact, through the calculation of Mathematica or Maple, we obtain that the number $a_{*}\approx 0.587459$ is the unique root in $(0, 1)$ of the equation $A_{2}^{\prime}(a)=0$, where
\begin{eqnarray*}
A_{2}^{\prime}(a)&=&2673-2502 a-2025 a^{2}+332 a^{3}+255 a^{4}-6 a^{5}-7 a^{6}\\
&&+162\times 14.796883 (a+1)^{3}(5 a-3).
\end{eqnarray*}
(In fact, the equation $A_{2}^{\prime}(a)=0$ has six real roots, and the other five roots are out of the interval $[0, 1]$).

Now we plug the value of $\lambda$ into the expression for  $A_{2}(a_*)$. This gives
$$A_{2}(a_*)=\frac{1}{5 a_*-3} \psi(a_*),
$$
where $\psi(a)$ is given by \eqref{liu21}. By Lemma \ref{HLP-lem5}, we know that $a_{*}$ is the unique root in $(0,1)$ of the equation $\psi(a)=0$. Thus, $A_{2}(a_{*})=0$ and $A_{2}^{\prime}(a_{*})=0$. Besides this observation, we have $A_{2}(0)\approx 2867.904954>0$ and $ A_{2}(1)=6144>0$. Consequently, $A_{2}(a)\geq 0$ in the interval $(0, 1)$, which proves that $A(r) \leq 1 $ for $r \leq 1/3$.

Finally, to prove that the constant $\lambda$ is sharp, we consider the function $f$ given by
\begin{eqnarray}
f(z)=\frac{a-z}{1-a z}=a-\left(1-a^{2}\right) \sum_{n=1}^{\infty} a^{n-1} z^{n},~z\in \ID,
\label{1}
\end{eqnarray}
where $a \in(0,1)$.

For this function, straightforward calculations show that
\begin{eqnarray*}
A_{\lambda _1}(r)&:=&\sum_{n=0}^{\infty}\left|a_{n}\right| r^{n}+\left(\frac{1}{1+\left|a_{0}\right|}+\frac{r}{1-r}\right) \sum_{n=1}^{\infty}\left|a_{n}\right|^{2} r^{2 n}+\frac{8}{9}\left(\frac{S_{r}}{\pi}\right)+\lambda_{1}\left(\frac{S_{r}}{\pi}\right)^{2}\\
&=& a+\frac{\left(1-a^{2}\right) r}{1-a r}+\frac{1+a r}{(1+a)(1-r)} \frac{\left(1-a^{2}\right)^{2} r^{2}}{1-a^{2} r^{2}}+\frac{8}{9} \frac{\left(1-a^{2}\right)^{2} r^{2}}{\left(1-a^{2} r^{2}\right)^{2}}+\lambda_{1} \frac{\left(1-a^{2}\right)^{4} r^{4}}{\left(1-a^{2} r^{2}\right)^{4}}\\
&=& a+\frac{\left(1-a^{2}\right) r}{1-r}+\frac{8}{9} \frac{\left(1-a^{2}\right)^{2} r^{2}}{\left(1-a^{2} r^{2}\right)^{2}}+\lambda_{1} \frac{\left(1-a^{2}\right)^{4} r^{4}}{\left(1-a^{2} r^{2}\right)^{4}}.
\end{eqnarray*}

For $r=1/3$, the last expression becomes
$$
A_{\lambda _1}(1/3)=a+\frac{1-a^{2}}{2}+\frac{8\left(1-a^{2}\right)^{2}}{\left(9-a^{2}\right)^{2}}+81 \lambda \frac{\left(1-a^{2}\right)^{4}}{\left(9-a^{2}\right)^{4}}+81\left(\lambda_{1}-\lambda\right) \frac{\left(1-a^{2}\right)^{4}}{\left(9-a^{2}\right)^{4}}.
$$
Choose $a$ as the positive root $a_{*}$ of the equation $\psi(a)= 0$. As a consequence, we see that
$$
A_{\lambda _1}(r)=1+81\left(\lambda_{1}-\lambda\right) \frac{\left(1-a^{2}\right)^{4}}{\left(9-a^{2}\right)^{4}}
$$
which is bigger than 1 in case $ \lambda_{1}>\lambda $. This proves the sharpness assertion and the proof of Theorem \ref{HLP-th1} is complete.\hfill $\Box$

\subsection{Proof of Theorem \ref{HLP-th2}}
Let $|a_{0}|=a \in(0,1) $. By Lemma \ref{HLP-lem1} with $N=1$, and  Lemma \ref{HLP-lem2}, we have
$$
B(r) \leq a^{2}+\frac{\left(1-a^{2}\right) r}{1-r}+\frac{9}{8} \frac{\left(1-a^{2}\right)^{2} r^{2}}{\left(1-a^{2} r^{2}\right)^{2}}+\mu  \frac{\left(1-a^{2}\right)^{4} r^{4}}{\left(1-a^{2} r^{2}\right)^{4}}:=B_{1}(r).
$$
Since $B_{1}(r)$ is an increasing function of $r$, we have for $r\leq 1/(3-a)$ that
\begin{eqnarray*}
B(r) &\leq & B_{1}\left(\frac{1}{3-a}\right)= a^{2}+\frac{1-a^{2}}{2-a}+\frac{(3-a)^{2}\left(1-a^{2}\right)^{2}}{8(3-2 a)^{2}}+\frac{\mu}{81} \frac{(3-a)^{4}\left(1-a^{2}\right)^{4}}{(3-2 a)^{4}}\\
&=&1-\frac{(1+a)(1-a)^{3}}{648(2-a)(3-2 a)^{4}} B_{2}(a),
\end{eqnarray*}
where
\begin{eqnarray*}
B_{2}(a)&=&39366-80919 a+59778 a^{2}-19197 a^{3}+2916 a^{4}-324 a^{5}\\
&&-8 \mu\left(162+27 a-378 a^{2}+30 a^{3}+290 a^{4}-108 a^{5}-62 a^{6}+50 a^{7}-12 a^{8}+a^{9}\right).
\end{eqnarray*}

Following the idea of the proof of Theorem \ref{HLP-th1}, we show that the function $B_{2}(a)$ has exactly one stationary point $a_{**}=0.638302\cdots$ in $[0, 1]$, which is the unique root in $[0, 1]$ of the equation $ B_{2}(a)=0$.

In fact, through the calculation of Mathematica or Maple, we obtain that the number $a_{**}=0.638302\cdots$ is the unique root in $(0, 1)$ of the equation $B_{2}^{\prime}(a)=0$, where
\begin{eqnarray*}
B_{2}^{\prime}(a)&=&80919-119556 a+57591 a^{2}-11664 a^{3}+1620 a^{4}\\
&&-8\times 13.966088 (a+1)^{2}(a-3)^{3}\left(9 a^{3}-33 a^{2}+29 a-1\right).
\end{eqnarray*}
(In fact, the equation $B_{2}^{\prime}(a)=0$ has eight real roots, and the other seven roots are out of the interval $[0, 1]$).



Now we plug the value of $\mu$ defined in the statement of the theorem into the  expression for $B_{2}(a_{**})$.
This gives
$$
B_{2}(a_{**})=\frac{1}{9 a_{**}^{3}-33 a_{**}^{2}+29 a_{**}-1} \phi(a_{**}).
$$

By Lemma \ref{HLP-lem6}, we know that $a_{**}$ is the unique root in $(0,1)$ of the equation $ \phi(a)=0$. Thus, $B_{2}(a_{**})=0$ and $B_{2}^{\prime}(a_{**})=0 $, by the choice of the value of $\mu$. Besides this observation, we have $B_{2}(0)\approx 21265.949952>0$ and $ B_{2}(1)=1620>0$. Consequently, $B_{2}(a)\geq 0$ in the interval $(0, 1)$, which proves that $B(r) \leq 1 $
for $r \leq 1/(3-a)$.

To prove that the constant $\lambda$ is sharp, we consider the function $f$ given by (\ref{1}) and compute
for this function the value of $B_{\mu _1}(r)$, which is obtained from $B(r)$ by replacing $\mu$ by $\mu_1$.
Indeed, straightforward calculations as before show that
\begin{eqnarray*}
B_{\mu _1}(r)= a^{2}+\frac{\left(1-a^{2}\right) r}{1-r}+\frac{9}{8} \frac{\left(1-a^{2}\right)^{2} r^{2}}{\left(1-a^{2} r^{2}\right)^{2}}+\mu_{1} \frac{\left(1-a^{2}\right)^{4} r^{4}}{\left(1-a^{2} r^{2}\right)^{4}}.
\end{eqnarray*}

For $r=1/(3-a)$, we have
\begin{eqnarray*}
B_{\mu _1}\left (\frac{1}{3-a}\right )&=&a^{2}+\frac{1-a^{2}}{2-a}+\frac{(3-a)^{2}\left(1-a^{2}\right)^{2}}{8(3-2 a)^{2}}+\frac{\mu}{81} \frac{(3-a)^{4}\left(1-a^{2}\right)^{4}}{(3-2 a)^{4}}\\
&&+\frac{1}{81}\left(\mu_{1}-\mu\right) \frac{(3-a)^{4}\left(1-a^{2}\right)^{4}}{(3-2 a)^{4}}.
\end{eqnarray*}
Choose $a$ as the unique root $a_{**}$ in $(0, 1)$ of the equation $\phi(a)= 0$. As a consequence, we find that
$$
B_{\mu _1}(r)=1+\frac{1}{81}\left(\mu_{1}-\mu\right) \frac{(3-a)^{4}\left(1-a^{2}\right)^{4}}{(3-2 a)^{4}}
$$
which is bigger than $1$ in case $ \mu_{1}>\mu $. Thus the proof of Theorem \ref{HLP-th2} is complete.\hfill $\Box$

\section{The Bohr-type operator on classes of subordination}\label{HLP-sec5}
\subsection{Proof of Theorem \ref{HLP-th3}}
Given $f\prec h$. Then $f(z)=h(\varphi(z))$ for some Schwarz function $\varphi$. We may let $h(z)=\sum_{n=0}^{\infty} b_{n} z^{n}$.
It follows that
$$f^{\prime}(z)=\varphi^{\prime}(z) h^{\prime}(\varphi(z))= \sum_{n=1}^{\infty} n b_{n} \varphi^{n-1}(z) \varphi^{\prime}(z)
~\mbox{ and }~s_{k}\left(f^{\prime}\right)=\sum_{n=1}^{k+1} n b_{n} s_{k}\left(\varphi^{n-1}(z) \varphi^{\prime}(z)\right).
$$

Now by \eqref{LHS-eq3}, we find that
$$\left|s_{k}\left(f^{\prime}\right)\right| \leq \sum_{n=1}^{k+1} n\left|b_{n}\right|\left|s_{k}\left(\varphi^{\prime}(z) \varphi^{n-1}(z)\right)\right| \leq \sum_{n=1}^{k+1} n\left|b_{n}\right||z|^{n-1}= B_{r}\left(s_{k}\left(h^{\prime}\right)\right)
$$
for $|z|\leq  \frac{1}{2}r_{0}(|\varphi '(0)|)$, which in particular proves the first part of the theorem, since
$r_{0}(|\varphi '(0)|) \geq 1-\sqrt{\frac{2}{3}}$.

For the second part, by the subadditivity of $ B_{r}$ and  \eqref{LHS-eq4}, we have
$$
 B_{r}\left(s_{k}\left(f^{\prime}\right)\right) \leq \sum_{n=1}^{k+1} n\left|b_{n}\right|  B_{r}\left(s_{k}\left(\varphi^{\prime} \varphi^{n-1}(z)\right)\right)\leq \sum_{n=1}^{k+1} n\left|b_{n} |\,| z\right|^{n-1}= B_{r}\left(s_{k}\left(h^{\prime}\right)\right)
$$
for $|z|\leq  r_{0}(|\varphi '(0)|)r_{1}(|\varphi '(0)|)$, which in particular proves the second part of the theorem, since
$r_{0}(|\varphi '(0)|)r_{1}(|\varphi '(0)|\geq \frac{1}{3}\Big(1-\sqrt{\frac{2}{3}}\Big)$. \hfill $\Box$
\vskip 2mm

%
%

\subsection{Proof of Theorem \ref{HLP-th6}}
Let $h(z)=\sum_{n=0}^{\infty} a_{n} z^{n}$. By Lemmas \ref{Theo-F} and \ref{HLP-lem3}, we have
\begin{eqnarray*}
 B_{r}\left((h \circ \varphi)^{\prime}\right) &\leq&  B_{r}\left(h^{\prime} \circ \varphi\right)  B_{r}\left(\varphi^{\prime}\right) \leq  B_{r}\left(h^{\prime}\circ \varphi\right) \leq
\sum_{n=1}^{\infty} n |a_{n} |\,|  B_{r}(\varphi^{n-1})|
\leq  B_{r}\left(h^{\prime}\right)
\end{eqnarray*}
for $0 \leq r \leq \frac{1}{3}\Big(1-\sqrt{\frac{2}{3}}\Big)$. \hfill $\Box$
\vskip 2mm

\subsection{Proof of Theorem \ref{HLP-th7}}
By assumption $f(z)=g(z) h(\varphi(z))$. Then
$$
f^{\prime}(z)=g^{\prime}(z) h(\varphi(z))+g(z) h^{\prime}(\varphi(z)) \varphi^{\prime}(z).
$$

By Theorems~F 
and \ref{HLP-th6}, we have
\begin{eqnarray*}
 B_{r}\left(f^{\prime}\right) &\leq&  B_{r}\left(g^{\prime}\right)  B_{r}(h\circ \varphi)+ B_{r}(g) B_{r}\left(h^{\prime}\circ \varphi\right)  B_{r}\left(\varphi^{\prime}\right)\\
&\leq&  B_{r}\left(g^{\prime}\right)  B_{r}(h)+ B_{r}(g)  B_{r}\left(h^{\prime}\right).
\end{eqnarray*}


Let $\phi(z) =\frac{g(z)-g(0)}{2b}$. Since $\phi(z)$ is an analytic self map on $\mathbb{D} $  and $\phi (0)=0$, Lemma \ref{Theo-F}
shows that
$$ B_{r}(g^{\prime})=2b B_{r}(\phi^{\prime})\leq 2b   ~\mbox{ for $|z|= r \leq r_{0}(|\varphi '(0)|)$}, ~r_0(x)=1-\sqrt{\frac{1+x}{2+x}}.
$$
Therefore,
\begin{eqnarray*}
 B_{r}\left(f^{\prime}\right) &\leq& b(2 B_{r}(h)+  B_{r}(h^{\prime})) ~\mbox{ for $|z|= r \leq r_{2}$},
\end{eqnarray*}
where $r_{2}=\min \left\{r_{1}(|g'(0)|), r_{0}(|\varphi '(0)|)\right\}\geq \min \left\{\frac{1}{3}, 1-\sqrt{\frac{2}{3}}\right\}= 1-\sqrt{\frac{2}{3}}$. \hfill $\Box$
\vskip 2mm

\subsection*{Acknowledgments}
This research of the first two authors are partly supported by Guangdong Natural Science Foundations (Grant No. 2021A1515010058).
The work of the third author was supported by Mathematical Research Impact Centric Support (MATRICS) of the Department of Science and Technology (DST), India  (MTR/2017/000367).

\subsection*{Conflict of Interests}
The authors declare that they have no conflict of interest, regarding the publication of this paper.

\subsection*{Data Availability Statement}
The authors declare that this research is purely theoretical and does not associate with any datas.


\begin{thebibliography}{14}

\bibitem{Abu} Y. Abu-Muhanna, Bohr's phenomenon in subordination and bounded harmonic classes,
\textit{Complex Var. Elliptic Equ.} \textbf{55}(11) (2010),  1071--1078.

\bibitem{Abu2} Y. Abu-Muhanna and R. M. Ali,
Bohr's phenomenon for analytic functions into the exterior of a compact convex body,
\textit{J. Math. Anal. Appl.} \textbf{379}(2) (2011), 512--517.

\bibitem{Abu4} Y. Abu-Muhanna and R. M. Ali,
Bohr's phenomenon for analytic functions and the hyperbolic metric,
\textit{Math. Nachr.} \textbf{286}(11-12) (2013), 1059--1065.

\bibitem{AAL20} Y. Abu-Muhanna, R. M. Ali and S. K. Lee,
Bohr operator on analytic functions, \textit{https://arxiv.org/abs/1912.11787 }

\bibitem{Abu3} Y. Abu-Muhanna, R. M. Ali, Z. C. Ng, and S. F. M. Hasni,
Bohr radius for subordinating families of analytic functions and bounded harmonic mappings,
\textit{J. Math. Anal. Appl.} \textbf{420}(1) (2014), 124--136.

\bibitem{AAP2016} Y. Abu-Muhanna,  R.~M. Ali,  and S. Ponnusamy, On the Bohr inequality, In Progress in Approximation Theory
and Applicable Complex Analysis, (Edited by N.K. Govil et al.),
\textit{Springer Optimization and Its Applications}, \textbf{117} (2016), 265--295.

\bibitem{A2012} L. Aizenberg, Remarks on the Bohr and Rogosinski phenomena for power series,
\textit{Anal. Math. Phys.} \textbf{2}(1) (2012), 69--78.

\bibitem{AES2005} L. Aizenberg, M. Elin and D. Shoikhet,
On the Rogosinski radius for holomorphic mappings and some of its aplications,
\textit{Studia Math.} \textbf{168}(2) (2005), 147--158.

\bibitem{AliBarSoly}  R. M. Ali, R. W. Barnard and A. Yu. Solynin,
A note on the Bohr's phenomenon for power series,
\textit{J. Math. Anal. Appl.} \textbf{449}(1) (2017), 154--167.

\bibitem{AlKayPON-19} S. A. Alkhaleefah, I. R. Kayumov and S. Ponnusamy,
On the Bohr inequality with a fixed zero coefficient,
\textit{Proc. Amer. Math. Soc.} \textbf{147}(12) (2019), 5263--5274.

\bibitem{BhowDas-18} B. Bhowmik and N. Das, Bohr phenomenon for subordinating families of certain univalent functions,
\textit{J. Math. Anal. Appl.} \textbf{462}(2) (2018), 1087--1098.

\bibitem{BD2019} B. Bhowmik and N. Das, A note on the Bohr inequality, {\tt https://arxiv.org/pdf/1911.06597.pdf}

\bibitem{BK1997} H.P. Boas and D. Khavinson, Bohr's power series theorem in several variables,
\textit{Proc. Amer. Math. Soc.} \textbf{125}(10) (1997), 2975--2979.

\bibitem{B1914}  H. Bohr, A theorem concerning power series, \textit{Proc. Lond. Math. Soc.} \textbf{2}(13) (1914), 1--5.



\bibitem{DFOOS} A. Defant, L. Frerick, J. Ortega-Cerd\`{a}, M. Ouna\"{i}es, and K. Seip,
The Bohnenblust-Hille inequality for homogenous polynomials is hypercontractive,
\textit{Ann. of Math. }  \textbf{174}(2) (2011), 512--517.


\bibitem{EPR-2017} S. Evdoridis, S. Ponnusamy and A. Rasila,
Improved Bohr's inequality for locally univalent harmonic mappings,
\textit{Indag. Math. (N.S.)}, \textbf{30} (2019), 201--213.


\bibitem{GarMasRoss-2018} S.~R.~Garcia, J.~ Mashreghi and W.~T.~Ross,
\textit{Finite Blaschke products and their connections}, Springer, Cham, 2018.

\bibitem{HLP2020} Y.~Huang, M.~S. Liu and S. Ponnusamy, Refined bohr-type inequalities with area measure for bounded analytic functions, \textit{Analysis and Mathematical Physics}, \textbf{10}(4) (2020), 21 pages, Aricle 50.

\bibitem{IKKP2018} A. A. Ismagilov, A. V. Kayumova, I. R. Kayumov and S. Ponnusamy,
Bohr inequalities in some classes of analytic functions. (Russian) Complex analysis (Russian), 69--83,
\textit{Itogi Nauki Tekh. Ser. Sovrem. Mat. Prilozh. Temat. Obz., 153, Vseross. Inst. Nauchn. i Tekhn. Inform. (VINITI), Moscow,} 2018.

\bibitem{IKKP2020} A. A. Ismagilov, I. R. Kayumov and S. Ponnusamy,
Sharp Bohr type  inequality, \textit{J. Math. Anal. and Appl.}, \textbf{486}(1)(2020), 10 pages; Article 124147.

\bibitem{KP2017}   I.~R. Kayumov and S. Ponnusamy, Bohr inequality for odd analytic functions,
\textit{Comput. Methods Funct. Theory}, \textbf{17} (4)(2017),  679--688.

\bibitem{KS2017}    I.~R. Kayumov and S. Ponnusamy, Bohr's inequality for analytic functions with lacunary series and harmonic functions,
\textit{J. Math. Anal. and Appl.}, 465(2) (2018),  857--871. 

\bibitem{KayPon3} I. R. Kayumov and S. Ponnusamy, Improved version of Bohr's inequality, \textit{C. R. Math. Acad. Sci. Paris}, \textbf{356}(3) (2018),  272--277.


\bibitem{KP2019} I. R. Kayumov and S. Ponnusamy, Bohr-Rogosinski radius for analytic functions,
{\tt  arXiv:1708.05585v1}

\bibitem{KSS2017}   I.~R. Kayumov, S. Ponnusamy, and N. Shakirov, Bohr radius for locally univalent harmonic mappings,
\textit{Math. Nachr}. \textbf{291} (2017),  1757--1768.

\bibitem{KKP2018} A.  Kayumova, I. R. Kayumov and S. Ponnusamy, Bohr's inequality for harmonic mappings and beyond,
\textit{Mathematics and computing}, 245--256, Commun. Comput. Inf. Sci., 834, Springer, Singapore, 2018.

\bibitem{EDB1986} E. Landau and D. Gaier, Darstellung und Begr\"{u}undung einiger neuerer Ergebnisse der Funktionentheorie,
Springer-Verlag, 1986.

\bibitem{LLP2020}  G. Liu, Z.H. Liu and S. Ponnusamy, Refined Bohr inequality for bunded analytic functions,
\textit{Bulletin des sciences math\'{e}matiques}, To appear; See also {\tt arXiv:2006.08930v1}



\bibitem{LP2018}  G. Liu and  S. Ponnusamy, {On Harmonic $\nu$-Bloch and $\nu$-Bloch-type mappings},
\textit{ Results Math.} \textbf{73}(3)(2018), Art. 90, 21 pages.

\bibitem{LP2020}  G. Liu and  S. Ponnusamy, Improved Bohr inequality for harmonic mappings,
\textit{Math. Nachr}. (2021), 14 pp. {\tt arXiv:submit/3671842}

\bibitem{LSX2018}  M.~S. Liu, Y.~M. Shang, and J. F. Xu, Bohr-type inequalities of analytic functions,
\textit{J. Inequal. Appl.}, \textbf{2018}(2018), Art.345, 13 pages.

\bibitem{LPW2020} M.~S. Liu, S. Ponnusamy and J. Wang, Bohr's phenomenon for the classes of Quasi-subordination and
$K$-quasiregular harmonic mappings, \textit{RACSAM} \textbf{114}(2020), Art. 115, 15 pages. https://doi.org/10.1007/s13398-020-00844-0

\bibitem{LP2019} Z. H. Liu and S. Ponnusamy, Bohr radius for subordination and $K$-quasiconformal harmonic mappings,
\textit{Bull. Malays. Math. Sci. Soc.}, \textbf{42} (2019), 2151--2168.


\bibitem{PPS2002}  V.~I. Paulsen, G. Popascu, and D. Singh, On Bohr's inequality,
\textit{Proc. Lond. Math.}, \textbf{85}(2) (2002),  493--512.

\bibitem{PS2004}  V.~I. Paulsen and D. Singh, Bohr's inequality for uniform algebras,
\textit{Proc. Amer. Math. Soc.}, \textbf{132} (2004),  3577--3579.

\bibitem{PS2006}   V.~I. Paulsen and D. Singh, Extensions of Bohr's inequality,
\textit{Bull. Lond. Math. Soc.}, \textbf{38}(6) (2006),  991--999.

\bibitem{PVW2019}  S. Ponnusamy, R. Vijayakumar and K.-J. Wirths,
New inequalities for the coefficients of unimodular bounded functions,
\emph{Results Math.} \textbf{75}, 107 (2020). {\tt https://doi.org/10.1007/s00025-020-01240-1}

\bibitem{PW2019}  S.~Ponnusamy and K.-J.~Wirths, Bohr type inequalities for functions with a multiple zero at the origin,
\textit{Comput. Methods Funct. Theory} \textbf{20}(2020), 559--570.



\bibitem{WR1923} W. Rogosinski, \"{U}ber Bildschranken bei Potenzreihen und ihren Abschnitten, \textit{Math. Z}, \textbf{17} (1923), 260--276.

\bibitem{IG1925} I. Schur und G. Szeg\"{o}, \"{U}die Abschnivse einer im Einheitskreise beschr\"{a}nkten Potenzreihe,
\textit{Sitz.-Ber.Preuss.Acad.Wiss.Berlin phys.-Math.Kl}, \textbf (1925), 545--560.

\bibitem{S1927} S. Sidon, \"{U}ber einen Satz von Herrn Bohr, \textit{Math. Z.} \textbf{26}(1) (1927), 731--732.

\bibitem{T1962}  M. Tomi\'c, Sur un th\'eor\`eme de H. Bohr, \textit{Math. Scand.} \textbf{11} (1962), 103--106.





\end{thebibliography}
\end{document}